\documentclass[12pt]{article}

\usepackage{amsmath,amsfonts,amssymb,amsthm,pstricks,pst-plot,pst-node,eucal}
\usepackage[dvips]{graphics}

\newcommand{\Z}{\mathbb{Z}}
\newcommand{\R}{\mathbb{R}}
\newcommand{\oh}{\tfrac12}
\newcommand{\fS}{\mathfrak{S}}
\newcommand{\tfS}{\widetilde{\mathfrak{S}}}
\newcommand{\fM}{\mathfrak{M}}
\newcommand{\bK}{\mathsf{K}}
\newcommand{\bE}{\mathsf{E}}
\newcommand{\bB}{\mathsf{B}}
\newcommand{\cL}{\mathcal{L}}
\newcommand{\bS}{\mathsf{S}}
\newcommand{\bfS}{\boldsymbol{\mathfrak{S}}}
\DeclareMathOperator{\Prob}{Prob} 
\DeclareMathOperator{\Ad}{Ad} 
\DeclareMathOperator{\Var}{Var} 
\DeclareMathOperator{\Ai}{Ai} 
\DeclareMathOperator{\Res}{Res} 

\newcommand{\cl}{\{\lambda(t)\}}
\newcommand{\vac}{v_\emptyset}
\newcommand{\Dt}{\Delta t}
\newcommand{\Dh}{\Delta h}
\newcommand{\Pl}{\text{Planch}}
\newcommand{\qd}[1]{(#1;q)_\infty}

\newcommand{\timeprod}{\operatornamewithlimits{\overleftarrow{\prod}}}

\newpsobject{showgrid}{psgrid}{subgriddiv=1,griddots=5,gridlabels=0pt}

\newtheorem{lemma}{Lemma}
\newtheorem{theorem}{Theorem}
\newtheorem{corollary}{Corollary}
\DeclareMathOperator{\dil}{dilog}

\newcommand{\LV}{\Lambda^{\frac\infty2}V}
\newcommand{\ul}{\underline}
\newcommand{\al}{\alpha}

\theoremstyle{definition}
\newtheorem{remark}{Remark}
\newtheorem{definition}{Definition}

\begin{document}

\title{Correlation function of Schur process with 
application to local geometry of  a random 
3-dimensional Young diagram}
\author{Andrei Okounkov and Nikolai Reshetikhin\thanks{
 Department of Mathematics, University of California at
Berkeley, Evans Hall \#3840, 
Berkeley, CA 94720-3840. E-mail: okounkov@math.berkeley.edu, 
reshetik@math.berkeley.edu }
}
\date{}
\maketitle 

\begin{abstract}
Schur process is a time-dependent analog of the Schur measure on
partitions studied in \cite{O}. Our first result is that the correlation 
functions of the Schur process are
determinants with a kernel that has a nice contour integral
representation in terms of the parameters of the process. 
This general result is then applied to a particular specialization
of the Schur process, namely to random 3-dimensional Young diagrams. 
The local geometry of a large random 3-dimensional diagram is described 
in terms of a determinantal point process on a 2-dimensional lattice with
 the incomplete beta function kernel (which generalizes the discrete sine 
kernel). A brief discussion of the universality of this answer
concludes the paper. 
\end{abstract}

\section{Introduction}

\subsection{Schur measure and Schur process} 

\subsubsection{}

This paper is a continuation of \cite{O}. The Schur measure, introduced
 in \cite{O}, is a measure on partitions $\lambda$ which weights a 
partitions $\lambda$ proportionally to $s_\lambda(X) \, s_\lambda(Y)$,
where $s_\lambda$ is the Schur function and $X$ and $Y$ are two 
sets of variables.

Here we consider
a time-dependent version of the Schur measure which we call
the Schur process, see Definition \ref{defS}. 
This is a measure on sequences $\cl$ such
that
$$
\Prob(\cl)\,\propto\,  \prod \bS^{(t)}(\lambda(t),
\lambda(t+1))\,,
$$ 
where the time-dependent weight $\bS^{(t)}(\lambda,\mu)$ is 
a certain generalization of a skew Schur function. It is given
by a suitably regularized infinite minor of a certain Toeplitz
matrix. This determinants can be also interpreted as a Karlin-McGregor
non-intersection probability \cite{KM}. The distribution of each 
individual $\lambda(t)$ is then a Schur measure with suitable
parameters.

In this paper, we consider the Schur process only in 
discrete time, although a continuous time formulation is 
also possible, see Section \ref{cont}.
 
The idea to make the Schur measure
time dependent was inspired, in part, by the paper 
\cite{J3} which deals with some particular instances of 
the general concept of the Schur processes introduced here.
For another development of the ideas of \cite{J3}, see
the paper \cite{PS} which appeared after the results of
the present paper were obtained.

\subsubsection{}
Our main interest in this paper are the correlation functions
of the Schur process which, by definition, are
the probabilities that the random set
$$
\fS( \{\lambda(t)\}) = \{(t,\lambda(t)_i-i+\oh)\}\,, \quad t,i\in Z\,,
$$
contains a given set $U\subset \Z \times (\Z+\oh)$. 

Using the same infinite wedge machinery as in \cite{O}
we prove that the correlation functions of the Schur
process have a determinantal form
$$
\Prob\left(U\subset \fS(\cl)\right) =
\det\big(K(u_i,u_j)\big)_{u_i,u_j\in U}\,,
$$
 and give an explicit
contour integral representation for the correlation kernel $K$,
see Theorem \ref{t1}. This theorem is a generalization
of Theorem 2 in \cite{O}. 

\subsection{Asymptotics}

\subsubsection{}

The contour integral representation for correlation
kernel is particularly convenient for
asymptotic investigations. Only very elementary
means, namely residue calculus and basic saddle-point
analysis are needed to derive the asymptotics. 

\subsubsection{}

We illustrate this in
Section \ref{secA} by computing explicitly the
correlation functions asymptotics for 3-dimensional
Young diagram in the bulk of their limit shape. 
Three-dimensional diagrams, also known as plane
partition, is an old subject which has received 
recently a lot of attention. 
More specifically, the
measure on 3D diagrams $\pi$ such that
\begin{equation}
\Prob(\pi)\, \propto \, q^{|\pi|}\,, \quad 0<q<1\,,
\label{qpi}
\end{equation}
where $|\pi|$ is the volume of $\pi$,  is 
a particular specialization of the Schur process for which 
the time parameter has the meaning of an extra spatial
dimension.
As $q\to 1$, a typical partition, suitably scaled,
 approaches a limit
shape which was described in \cite{CK} and can be seen 
in Figure \ref{f7} below. The existence of this limit
shape and some of its properties were
first established by A.~Vershik \cite{V},
who used direct combinatorial methods.

It is well known
that 3D diagrams are in bijection with certain rhombi
tiling of the plane, see for example Section 
\ref{rtil} below.  Local statistics of random domino tilings
have been the subject of intense recent studies, see
for example \cite{BP,CEP,CKP,K1} and the survey \cite{K2}.

The authors of \cite{CKP} 
computed local correlations for domino tilings
with periodic boundary conditions
and conjectured that the same formula holds in the 
thermodynamic limit for domino tilings of more 
general regions, see Conjecture 13.5 in \cite{CKP}. 

{}From the point of view of statistical mechanics,
the argument 
behind this conjecture can be the belief that 
in the thermodynamic limit and away from the 
boundary the local correlations depend 
only on macroscopic parameters, 
such as the density of tiles of a given kind.
In particular, in the case of 3D diagrams the densities
of rhombi of each of the 3 possible kinds
are uniquely fixed by  the 
tilt of the limit shape at the point in question.

\subsubsection{}

In this paper, we approach the local geometry of
random 3D diagrams using the general exact formulas
for correlation functions of the Schur process. For 3D diagrams, these 
general formulas specialize to contour integrals 
involving the quantum dilogarithm function \eqref{qdil}, see Corollary 
\ref{cor3D}.

We compute the $q\to 1$ limit of the
correlation kernel explicitly. The result turns out to 
be  the discrete incomplete beta kernel,
see Theorem \ref{t2}. One can show, see Section \ref{iint}, that this kernel is
 a specialization of the kernel of \cite{CKP} when one of the 
parameters vanishes in 
agreement with Conjecture 13.5 of \cite{CKP}.

The incomplete beta kernel is also a bivariate generalization of the 
discrete sine kernel which appeared in \cite{BOO} in the 
situation of the Plancherel specialization of the 
Schur measure, see also \cite{J3}.

\subsection{Universality}

Although we focus on 
one specific asymptotic problem, our methods are both
very basic and completely general. 
They should apply, therefore, with little or no 
modification in a much wider variety of situations
yielding the same or analogous results. 
In other words, both the methods and the results
should be universal in a large class of asymptotic
problems. 

As explained in Section \ref{eqtime},
this universality should be especially robust for 
the equal time correlations, in which case the 
discrete sine kernel should appear.
This is parallel to the situation with random 
matrices \cite{J2}. 

\subsection{Acknowledgments}

We are grateful to R.~Kenyon and A.~Vershik for fruitful
discussions. A.O.\ was partially supported by NSF grant DMS-0096246
and a Sloan foundation fellowship. N.~R.\ 
was partially supported by NSF grant DMS-0070931.
Both authors were partially supported by CRDF grant RM1-2244.

\section{Schur process}

\subsection{Configurations}

\subsubsection{} 

Recall that a partition is a sequence
$$
\lambda=(\lambda_1 \ge \lambda_2 \ge \lambda_3 \ge \dots \ge 0)
$$
of integers such that $\lambda_i=0$ for $i\gg 0$. The zero, or
empty, partition is denoted by $\emptyset$. The book \cite{M}
is a most comprehensive reference on partitions and symmetric functions. 

Schur process is a measure on sequences 
$$
\{ \lambda(t)\}\,, \quad t \in \Z\,,
$$
where each $\lambda(t)$ is a partition and $\lambda(t)=\emptyset$
for $|t|\gg 0$. We call the variable $t$ the time, even though 
it may have a different interpretation in applications. 

\subsubsection{} 

\begin{figure}[htbp]
  \begin{center}
    \scalebox{0.3}{\includegraphics{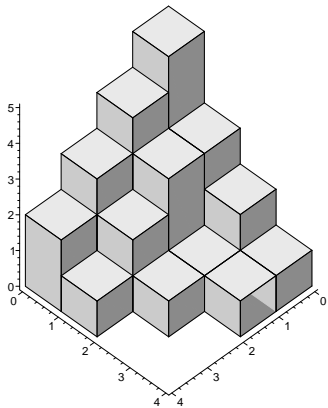}} 
    \caption{A 3-dimensional diagram $\pi$}
    \label{f1}
  \end{center}
\end{figure}

An example of such an object is a plane partition
which, by definition, is a 2-dimensional array of nonnegative numbers 
$$
\pi=(\pi_{ij})\,, \quad i,j=1,2,\dots\,,
$$
that are nonincreasing as function of both $i$ and $j$ and such that
$$
|\pi|=\sum \pi_{ij} 
$$
is finite. The plot of the function
$$
(x,y)\mapsto \pi_{\lceil x \rceil,\lceil y \rceil}\,, \quad x,y>0\,,
$$
is a 3-dimensional Young diagram with volume $|\pi|$. 
For example, Figure \ref{f1} shows the 
3D diagram corresponding to the plane partition
\begin{equation}
  \label{pi}
 \pi = 
\left (\begin {array}{cccc} 5&3&2&1\\\noalign{\medskip}4&3&1&1
\\\noalign{\medskip}3&2&1\\\noalign{\medskip}2&1\end {array}
\right )\,,
 \end{equation}
where the entries that are not shown are zero.

We associate to $\pi$ the sequence $\cl$ of its diagonal slices, 
that is, the sequence of partitions
\begin{equation}
  \label{pila}
 \lambda(t)=(\pi_{i,t+i})\,, \quad i\ge\max(0,-t) \,. 
\end{equation}
It is easy to see that a configuration $\cl$ corresponds
to a plane partition if and only if it satisfies the 
conditions 
\begin{equation}
  \label{laint}
 \dots \prec \lambda(-2) \prec \lambda(-1) \prec \lambda(0)
\succ \lambda(1) \succ \lambda(2) \succ \dots\,, 
\end{equation}
where $\lambda \succ \mu$ means that $\lambda$ and $\mu$ interlace,
that is, 
$$
\lambda_1 \ge \mu_1 \ge \lambda_2 \ge \mu_2 \ge \lambda_3 \ge \dots \,.
$$
In particular, the configuration $\cl$ corresponding to 
the diagram \eqref{pi} is 
$$
(2) \prec (3,1) \prec (4,2) \prec (5,3,1) \succ (3,1)\succ (2,1)\succ (1)\,.
$$

\subsubsection{} 

The mapping 
$$
\lambda \mapsto \fS(\lambda)=\{\lambda_i - i +\oh\} \subset \Z+\oh
$$
is a bijection of the of the set of partitions and the 
set $\bfS$ of  subsets $\fS\subset \Z+\oh$ such that
$$
| \fS \setminus (\Z+\oh)_{<0} | = | (\Z+\oh)_{<0} \setminus \fS | < 
\infty \,.
$$
The mapping
\begin{equation}
  \label{fS}
  \{\lambda(t)\}\mapsto
\fS( \{\lambda(t)\}) = \{(t,\lambda(t)_i-i+\oh)\} \subset \Z \times (\Z+\oh)\,,
\end{equation}
identifies the configurations of the Schur process 
with certain subsets of $\Z \times (\Z+\oh)$. In other words, the mapping
\eqref{fS} makes the Schur process a  random point field on
$\Z \times (\Z+\oh)$. 

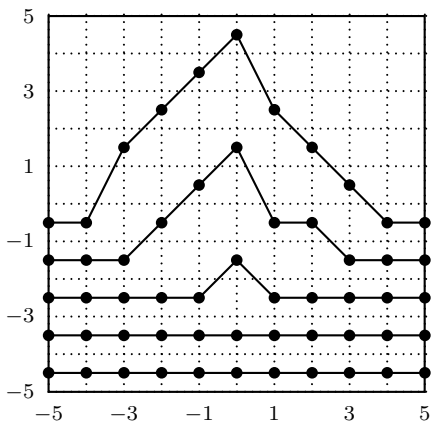
\begin{figure}[htbp]\psset{unit=0.5cm}
  \begin{center}
    \begin{pspicture}(-5,-5)(5,5)\showgrid
\scriptsize
\psaxes[axesstyle=frame,Ox=-5,Oy=-5,Dx=2,Dy=2,ticks=none](-5,-5)(5,5)
\psdots[dotscale=1.5](-5,-.5)(-5,-1.5)(-5,-2.5)(-5,-3.5)(-5,-4.5)(-4,-.5)(-4,-1.5)(-4,-2.5)(-4,-3.5)(-4,-4.5)(-3,1.5)(-3,-1.5)(-3,-2.5)(-3,-3.5)(-3,-4.5)(-2,2.5)(-2,-.5)(-2,-2.5)(-2,-3.5)(-2,-4.5)(-1,3.5)(-1, .5)(-1,-2.5)(-1,-3.5)(-1,-4.5)(0,4.5)(0,1.5)(0,-1.5)(0,-3.5)(0,-4.5)(1,2.5)(1,-.5)(1,-2.5)(1,-3.5)(1,-4.5)(2,1.5)(2,-.5)(2,-2.5)(2,-3.5)(2,-4.5)(3, .5)(3,-1.5)(3,-2.5)(3,-3.5)(3,-4.5)(4,-.5)(4,-1.5)(4,-2.5)(4,-3.5)(4,-4.5)(5,-.5)(5,-1.5)(5,-2.5)(5,-3.5)(5,-4.5)

\psline(-5,-.5)(-4,-.5)(-4,-.5)(-3,1.5)(-3,1.5)(-2,2.5)(-2,2.5)(-1,3.5)(-1,3.5)(0,4.5)(0,4.5)(1,2.5)(1,2.5)(2,1.5)(2,1.5)(3, .5)(3, .5)(4,-.5)(4,-.5)(5,-.5)
\psline(-5,-1.5)(-4,-1.5)(-4,-1.5)(-3,-1.5)(-3,-1.5)(-2,-.5)(-2,-.5)(-1, .5)(-1, .5)(0,1.5)(0,1.5)(1,-.5)(1,-.5)(2,-.5)(2,-.5)(3,-1.5)(3,-1.5)(4,-1.5)(4,-1.5)(5,-1.5)
\psline(-5,-2.5)(-4,-2.5)(-4,-2.5)(-3,-2.5)(-3,-2.5)(-2,-2.5)(-2,-2.5)(-1,-2.5)(-1,-2.5)(0,-1.5)(0,-1.5)(1,-2.5)(1,-2.5)(2,-2.5)(2,-2.5)(3,-2.5)(3,-2.5)(4,-2.5)(4,-2.5)(5,-2.5)
\psline(-5,-3.5)(-4,-3.5)(-4,-3.5)(-3,-3.5)(-3,-3.5)(-2,-3.5)(-2,-3.5)(-1,-3.5)(-1,-3.5)(0,-3.5)(0,-3.5)(1,-3.5)(1,-3.5)(2,-3.5)(2,-3.5)(3,-3.5)(3,-3.5)(4,-3.5)(4,-3.5)(5,-3.5)
\psline(-5,-4.5)(-4,-4.5)(-4,-4.5)(-3,-4.5)(-3,-4.5)(-2,-4.5)(-2,-4.5)(-1,-4.5)(-1,-4.5)(0,-4.5)(0,-4.5)(1,-4.5)(1,-4.5)(2,-4.5)(2,-4.5)(3,-4.5)(3,-4.5)(4,-4.5)(4,-4.5)(5,-4.5)
\end{pspicture}
    \caption{Point field and nonintersecting paths corresponding to the diagram in Figure \ref{f1}}
    \label{f2}
  \end{center}
\end{figure}

For example, the subset $\fS( \{\lambda(t)\})$
corresponding to the 3D diagram from Figure \ref{f1} is shown in 
Figure \ref{f2}. One can also visualize $\fS( \{\lambda(t)\})$ as a
collection of nonintersecting paths as in Figure \ref{f2}. 

\subsection{Probabilities}

\subsubsection{} 

Schur process, the formal definition of which is given 
in Definition \ref{defS} below, is a measure on sequences $\cl$ such that 
$$
\Prob(\{\lambda(t)\}) \,\propto\,
\prod_{t\in \Z} \bS^{(t)}(\lambda(t),\lambda(t+1))\,,
$$
where $\bS^{(t)}(\mu,\lambda)$ is a certain time-dependent transition 
weight between the partition $\mu$ and $\lambda$ which will be
defined presently. 

The coefficient $\bS^{(t)}(\mu,\lambda)$
 is a suitably regularized infinite minor
of a certain  Toeplitz matrix. It can be viewed as a generalization of the 
Jacobi-Trudy determinant for a skew Schur function or as a form 
of Karlin-McGregor non-intersection probability \cite{KM}. 

\subsubsection{}

Let a function 
$$
\phi(z) = \sum_{k\in\Z} \phi_k \, z^k\,,
$$
be nonvanishing on the unit circle
$|z|=1$ with winding number $0$ and geometric mean $1$.
These two conditions mean that $\log \phi$ is a well defined
function with mean $0$ on the unit circle. For simplicity,
we additionally assume that $\phi$
is analytic in some neighborhood of the
unit circle. Many of the results below
hold under weaker assumptions on $\phi$ 
as can be seen, for example, by an approximation
argument. We will not pursue here the greatest analytic
generality.

\subsubsection{}

Given two subsets 
$$
X=\{x_i\},\, Y=\{y_i\} \in \bfS\,,
$$
we wish to assign a meaning to the following infinite
determinant:
\begin{equation}
  \label{dT}
\det(\phi_{y_i-x_j}) \,.  
\end{equation}
We will see that, even though there is no canonical way
to evaluate this determinant, different regularizations
differ by a constant which depends only on 
$\phi$ and not on $X$ and $Y$.

Since our goal is to 
define probabilities only up to a constant factor, 
it is clear that different regularizations lead to 
the same random process. 

\subsubsection{} 

The function $\phi$ admits a Wiener-Hopf factorization 
of the form 
$$
\phi(z)=\phi^+(z)\,  \phi^-(z) 
$$
where the functions
$$
\phi^\pm(z)= 1 + \sum_{k\in\pm\mathbb{N}}
 \phi^\pm_{k} \, z^{k} 
$$
are analytic and nonvanishing in some neighborhood of the interior (resp., 
exterior) of the unit disk. 

There is a special case when the meaning of \eqref{dT} in 
unambiguous, namely, if $\phi^-(z)=1$ (or $\phi^+(z)=1$)
then the matrix is almost unitriangular and the determinant
in \eqref{dT} is essentially a finite determinant. This 
determinant is then a Jacobi-Trudy determinant for a skew Schur 
function
$$
\det\left(\phi^+_{\lambda_i-\mu_j+j-i}\right) = s_{\lambda/\mu} (\phi^+)\,,
$$
where $s_{\lambda/\mu} (\phi^+)$ is the skew Schur function 
$s_{\lambda/\mu}$ specialized so that
\begin{equation}
  \label{phih}
 h_k = \phi^+_k\,, 
\end{equation}
where $h_k$ are the complete homogeneous symmetric functions. 
Note that
\begin{equation}
  \label{svan}
  s_{\lambda/\mu}(\phi^+)=0\,, \quad \mu\not\subset\lambda \,.
\end{equation}

\subsubsection{}

\begin{definition}\label{defW}
  We define the transition weight by the following formula
\begin{equation}
  \label{T}
\bS_\phi(\mu,\lambda)=\sum_\nu s_{\mu/\nu}(\phi^-)\, s_{\lambda/\nu}(\phi^+) \,,
\end{equation}
where $s_{\lambda/\nu}(\phi^+)$ is defined by \eqref{phih} and $s_{\mu/\nu} (\phi^-)$ is the skew Schur function $s_{\mu/\nu}$
specialized so that
\begin{equation*}
 h_k = \phi^-_{-k}\,.  
\end{equation*}
\end{definition}

Note that if the determinant $\det(\phi_{\lambda_i-\mu_j+j-i})$ were
unambiguously defined and satisfied the Cauchy-Binet formula then 
it would equal \eqref{T} because 
$$
\big(\,\phi_{i-j}\,\big) = \big(\,\phi^+_{i-j}\,\big) \,  \big(\,\phi^-_{i-j}\,\big)\,.
$$
Also note that because of \eqref{svan}
the sum in \eqref{T} is finite, namely, it ranges over all 
$\nu$ such that $\nu\subset\mu,\lambda$.  

\subsubsection{}

In what follows, we will assume that the reader is familiar
with the basics of the infinite wedge formalism, see for 
example Chapter 14 of the book \cite{K} by V.~Kac.
An introductory account of this formalism, together
with some probablistic applications, can be found
in the lectures \cite{O2}.  We will
use the notation conventions of the appendix to \cite{O}
summarized for the reader's convenience in the appendix to
this paper. 

Consider the following vertex operators
$$
\Gamma_\pm(\phi)=\exp\left(
\sum_{k=1}^\infty (\log \phi)_{\mp k} \, \, \alpha_{\pm k} \right)\,,
$$
where $(\log \phi)_k$ denotes the coefficient of $z^{k}$
in the Laurent expansion of $\log \phi (z)$. In particular, if 
the algebra of symmetric functions is specialized as in 
\eqref{phih} then 
$$
(\log \phi)_k = \frac{p_k}{k}\,,\quad k=1,2,3,\dots\,,
$$
where $p_k$ is the $k$th power-sum symmetric function.

It is well known, see for example Exercise 14.26 in \cite{K},
that the matrix coefficient of the operators $\Gamma_\pm$
are the skew Schur functions, namely
$$
(\Gamma_-(\phi) \, v_\mu,
v_\lambda) = s_{\lambda/\mu}(\phi^+)\,.
$$
It follows that 
\begin{equation}
  \label{Tver}
  \bS_\phi(\mu,\lambda) = (\Gamma_-(\phi) \, \Gamma_+(\phi) \, v_\mu,
v_\lambda) \,.
\end{equation}
It is clear from the commutation relation 
\begin{equation}
  \label{commG}
  \Gamma_+(\phi)  \Gamma_-(\phi) = 
e^{\sum  k \, (\log \phi)_k (\log \phi)_{-k}} \, 
\Gamma_-(\phi) \, \Gamma_+(\phi) 
\end{equation}
that different ordering prescriptions in \eqref{Tver},
which produce different regularizations of \eqref{dT}, all
differ by a constant independent of $\mu$ and $\lambda$.

\subsubsection{}

Now suppose that for any half-integer $m\in \Z+ \oh$ we choose,
independently, a function 
$$
\phi[m](z)=\sum_{k\in \Z} \phi_k[m] \, z^k 
$$
as above so that the series
$$
\sum_{m\in\Z} \log\phi[m](z)
$$
converges absolutely and uniformly in some neighborhood of the unit disk.
This assumption is convenient but can be weakened.
 The functions  $\phi[m]$
will be the parameters of the Schur process. 

\begin{definition}\label{defS}
The
probabilities of the Schur process are given by
$$
\Prob(\{\lambda(t)\}) = \frac1Z \, \prod_{m\in \Z+1/2} 
\bS_{\phi[m]}\left(\lambda(m-\oh),\lambda(m+\oh)\right)\,,
$$
where the transition weight $\bS_\phi$ is defined in 
Definition \ref{defW} and $Z$ is the normalizing factor (partition function)
$$
Z=\sum_{\cl} \prod_{m\in \Z+1/2} 
\bS_{\phi[m]}\left(\lambda(m-\oh),\lambda(m+\oh)\right) \,.
$$
\end{definition}

\subsubsection{}

It follows from \eqref{Tver} that $Z$ is given by the 
following matrix coefficient
\begin{equation}
  \label{Zver}
  Z=\left(\timeprod_{m \in \Z+1/2} 
\Gamma_-(\phi[m]) \, \Gamma_+(\phi[m]) 
\, v_\emptyset, v_\emptyset\right) \,,
\end{equation}
where $\timeprod$ denotes the time-ordered product, that is,
the product in which operators are ordered from right to 
left in the increasing time order.

Using \eqref{commG} and the following consequence of 
\eqref{ann}
\begin{equation}
  \label{fix}
  \Gamma_+(\phi) \, \vac = \vac \,.
\end{equation}
we compute the matrix  coefficient \eqref{Zver}
as follows
$$
Z=\exp\left(\sum_{m_1<m_2} \sum_k k \, (\log\phi[m_1])_k \,
(\log\phi[m_2])_{-k} \right)\,.
$$
Our growth assumptions on the functions $\phi[m]$ 
ensure the convergence of $Z$.

\subsubsection{} 
It is clear from the vertex-operator description that
\begin{equation}
\Prob(\lambda(t)=\mu) \, \propto \, s_\mu\left(\prod_{m<t} \phi^+[m]
\right) \, s_\mu\left(\prod_{m>t} \phi^-[m]
\right) \,,\label{Schr}
\end{equation}
where, for example, the first factor is the image of the 
Schur function under the specialization that sets $h_k$ to 
the coefficient of $z^k$ in the product of $\phi^+[m](z)$
over $m<t$. 

This means that the distribution of each individual $\lambda(t)$
is a Schur measure \cite{O} with parameters \eqref{Schr}. 

\subsubsection{}\label{restrict}

More generally, the restriction of the Schur process to 
any subset of times is again a Schur process with suitably
modified parameters. Specifically, let a subset
$$
\{t_k\} \subset \Z\,,\quad k\in \Z\,,
$$
be given and consider the restriction of the Schur process to this
set, that is, consider the process
$$
\{\widetilde\lambda(k)\}=\{\lambda(t_k)\}\,.
$$
It follows from the vertex operator description that
this is again a Schur process with parameters
$$
\widetilde\phi[l]=\prod_{t_{l-\frac12} < m <t_{l+\frac12}} \phi[m]\,,
\quad l,m\in\Z+\oh\,.
$$
The case of a finite set $\{t_k\}$ is completely analogous
and can be dealt with formally by allowing infinite values
of $t_k$'s. 

\subsubsection{}\label{cont}

The restriction property from Section \ref{restrict}
forms a natural basis for
considering the Schur process in continuous time. 

For this we need a function $\cL(z,s)$ of a 
continuous variable $s\in\R$ which will 
play the role of the density of $\log\phi$. The restriction
$\{\lambda(t_k)\}$ 
of the process $\lambda(t)$, $t\in\R$, to any discrete 
set of times $\{t_k\}$ is the Schur process with 
parameters
$$
\phi[m](z) = \exp\left(\int_{t_{m-\frac12}}^{t_{m+\frac12}} \cL(z,s) \, ds\right)\,,
\quad m\in \Z+\oh\,.
$$

\subsubsection{}

Let $q\in (0,1)$ and consider the probability measure 
$\fM_q$ on the 
set of all 3D diagrams such that 
$$
\Prob(\pi) \,\propto\, q^{|\pi|} \,.
$$
We claim that there exists a particular choice of the 
parameters of the Schur process which yield this measure
under the correspondence \eqref{pila}. Concretely, set
\begin{equation}
  \label{specp}
\phi_{3D}[m](z)=
\begin{cases}
(1-q^{|m|} z)^{-1}\,,  & m<0\,,\\
(1-q^{|m|} z^{-1})^{-1}\,, & m>0 \,,
\end{cases}
\qquad m\in\Z+\oh\,.  
\end{equation}
It is well known that if the algebra of the symmetric
functions is specialized so that
$$
h_k = c^k \,, \quad k=1,2,\dots \,,
$$
for some constant $c$, then 
$$
s_{\lambda/\mu}=
\begin{cases}
c^{|\lambda|-|\mu|}\,, &  \mu \prec \lambda\,,\\
0\,, & \mu \not\prec \lambda\,.
\end{cases}
$$
Therefore, we have
\begin{multline*}
\bS_{\phi_{3D}[m]}\left(\lambda(m-\oh),\lambda(m+\oh)\right)=\\
  \begin{cases}
    q^{m(|\lambda(m-\frac12)|-|\lambda(m+\frac12)|)}\,, &
m<0,\, \lambda(m-\frac12) \prec \lambda(m+\frac12)\quad\textup{or}\\
 &
m>0,\, \lambda(m-\frac12) \succ \lambda(m+\frac12)\,,\\
0\,, & \textup{otherwise} \,.
  \end{cases}  
\end{multline*}
It follows that for the specialization 
\eqref{specp} the Schur process is
supported on configurations of the shape \eqref{laint} and
the weight of a configuration \eqref{laint} is proportional to 
$$
q^{\sum  |\lambda(t)|}=q^{|\pi|} \,.
$$
In particular, the partition function $Z$ becomes
\begin{multline*}
 Z_{3D}=\exp\left(\sum_{m_1,m_2=1/2}^\infty \sum_k \frac{q^{k(m_1+m_2)}}k \,\right)=\\
\prod_{m_1,m_2} (1-q^{m_1+m_2})^{-1} = \prod_{n=1}^\infty (1-q^n)^{-n} \,,
\end{multline*}
which is the well-known generating function, due to McMahon,
for 3D diagrams.

\subsection{Correlation functions}

\subsubsection{}

\begin{definition}
 Given a subset $U\subset \Z\times (\Z+\oh)$, define
the corresponding correlation function by
$$
\rho(U)=\Prob\left(U\subset \fS(\cl)\right)\,.
$$ 
These correlation functions depend on the parameters $\phi[m]$ of the Schur process. 
\end{definition}

In this section we show that 
\begin{equation}
  \label{rhodet}
  \rho(U)=\det\big(K(u_i,u_j)\big)_{u_i,u_j\in U}\,,
\end{equation}
for a certain kernel $K$ which will be computed explicitly. 

\subsubsection{}

Suppose that
$$
U=\{u_1,\dots,u_n\}\,, \quad u_i=(t_i,x_i)\in \Z\times (\Z+\oh)\,,
$$
and the points $u_i$ are ordered so that
$$
t_1 \le t_2 \le t_3 \dots \le t_n \,.
$$
For convenience, we set
$t_0 =- \infty$.  From \eqref{Tver} 
and \eqref{yn} it is clear that
\begin{equation}
  \label{rver}
  \rho(U)=\frac1Z\, (R_U \, v_\emptyset,v_\emptyset)\,,
\end{equation}
where $R_U$ is the following operator
$$
R_U = \timeprod_{m>t_n} \Gamma_-(\phi[m]) \, \Gamma_+(\phi[m]) 
\timeprod_{i=1..n} \left(\psi_{x_i} \psi_{x_i}^* \, 
\timeprod_{t_{i-1} < m < t_i} \Gamma_-(\phi[m]) \, \Gamma_+(\phi[m]) 
\right)
$$

\subsubsection{}

Define the operator
$$
\Psi_{x}(t)=\Ad\left(\prod_{m>t} \Gamma_+(\phi[m]) \,
\prod_{m<t} \Gamma_-(\phi[m])^{-1} \right) \cdot \psi_x \,,
$$
where $\Ad$ denotes the action by conjugation, and define
the operator $\Psi^*_x(t)$ similarly. Note the ordering of
the vertex operators inside the $\Ad$ symbol is immaterial
because the vertex operators commute up to a central element.

It follows from \eqref{rver}, \eqref{fix}, and \eqref{Zver} that
$$
\rho(U)=\left(\timeprod \Psi_{x_i}(t_i)\, \Psi^*_{x_i}(t_i) \, \vac,\vac
\right) \,.
$$
Now we apply Wick formula in the following form

\begin{lemma}[Wick formula]
Let $A_i = \sum_k a_{i,k} \, \psi_k$ and 
$A^*_i = \sum_k a^*_{i,k}\, \psi^*_k$. Then
\begin{equation}
  \label{Wick}
  \left(\prod A_i \, A^*_i \, \vac, \vac\right)=
\det\big(K_A(i,j)\big)\,,
\end{equation}
where
$$
K_A(i,j)=
\begin{cases}
(A_i \, A^*_j \, \vac,\vac)\,, & i \ge j\,,\\
-(A^*_j \, A_i \, \vac, \vac)\,, & i < j \,.
\end{cases}
$$
\end{lemma}

\begin{proof}
 Both sides of \eqref{Wick} are linear in $a_{i,k}$ and
$a^*_{i,k}$, therefore it suffices to verify
\eqref{Wick} for some linear basis in the space of possible
$A_i$'s and $A^*_j$'s. A convenient linear basis is formed
by the series \eqref{psiz} as the parameter $z$ varies. Using the canonical
anticommutation relation satisfied by $\psi(z)$ and $\psi^*(z)$ one then 
verifies \eqref{Wick} directly. 
\end{proof}

We obtain the formula \eqref{rhodet} with
$$
K((t_1,x_1),(t_2,x_2))=
\begin{cases}
  \left(\Psi_{x_1}(t_1)\, \Psi^*_{x_2}(t_2) \, \vac,\vac\right)\,,
& t_1 \ge t_2 \,,\\
-\left(\Psi^*_{x_2}(t_2)\, \Psi_{x_1}(t_1) \, \vac,\vac\right) \,,
& t_1 < t_2 \,.\\
\end{cases}
$$
Note that for $t_1\ne t_2$ the operators 
$\Psi_{x_1}(t_1)$ and $\Psi^*_{x_2}(t_2)$ do not, in general, 
anticommute so the time ordering is important. However,
for $t_1=t_2$ and $x_1\ne x_2$, these operators do 
anticommute, so at equal time the ordering is immaterial.

\subsubsection{}

A convenient generating function for the kernel $K$ can be 
obtained as follows. Set
\begin{equation}
\psi(z)=\sum_{k\in\Z+1/2} z^k \, \psi_k \,, 
\quad 
\psi^*(z)=\sum_{k\in\Z+1/2} z^{-k} \, \psi^*_k \,.\label{psiz}
\end{equation}
Set also 
$$
\Psi(t,z)=\Ad\left(\prod_{m>t} \Gamma_+(\phi[m]) \,
\prod_{m<t} \Gamma_-(\phi[m])^{-1} \right) \cdot \psi(z) \,,
$$
and define $\Psi^*(t,z)$ similarly. 

We have from \eqref{adal}
\begin{alignat*}{2}
  \Ad(\Gamma_\pm(\phi))\,\cdot\, &\psi(z) &&= \phi^\mp(z^{-1})\, \psi(z)\,,\\
\Ad(\Gamma_\pm(\phi))\,\cdot\, &\psi^*(z) &&= \phi^\mp(z^{-1})^{-1}\, \psi(z)\,.
\end{alignat*}
and therefore
\begin{equation*}
  \Psi(t,z)=\Phi(t,z)  \,
\psi(z)\,,  \quad
\Psi^*(t,z)=\Phi(t,z)^{-1} \,
\psi(z)\,, 
\end{equation*}
where
\begin{equation}
  \label{Phi}
 \Phi(t,z)=\frac{\prod_{m>t} \phi^-[m](z^{-1})}{\prod_{m<t} \phi^+[m](z^{-1})} \,. 
\end{equation}
Finally, it is obvious from the definitions that
\begin{align*}
 (\psi(z) \, \psi^*(w) \, \vac,\vac)&=
\frac{\sqrt{zw}}{z-w}\,, \quad |z|>|w| \,, \\
 -(\psi^*(w) \, \psi(z) \, \vac,\vac)&=
\frac{\sqrt{zw}}{z-w}\,, \quad |z|<|w| \,. 
\end{align*}
Putting it all together, we obtain the following

\begin{theorem}\label{t1}
We have
$$
 \rho(U)=\det\big(K(u_i,u_j)\big)_{u_i,u_j\in U}\,,
$$
where the kernel $K$ is determined by the following
generating function
\begin{align}
\bK_{t_1,t_2}(z,w)&=\sum_{x_1,x_2\in\Z+\frac12} z^{x_1} \, w^{-x_2} \,
K((t_1,x_1),(t_2,x_2))\,, \label{bK1}\\
&=\frac{\sqrt{zw}}{z-w} \, \frac{\Phi(t_1,z)}
{\Phi(t_2,w)}\,. \label{bK2}
\end{align}
Here the function $\Phi(t,z)$ is defined by \eqref{Phi}, and
\eqref{bK1} is the expansion of \eqref{bK2}
in the region $|z|>|w|$ if  $t_1 \ge t_2$ and
$|z|<|w|$ if $t_1 < t_2$. 
\end{theorem}

\subsubsection{}

In the special case of the measure $\fM_q$ on the 
3D diagrams the function 
\eqref{Phi} specializes to the following function
\begin{equation}
  \label{Phip}
  \Phi_{3D}(t,z)=\frac
{\prod_{m > \max(0,-t)} (1-q^{m}/z)}
{\prod_{m > \max(0,t)} (1-q^{m} z)}\,, \quad m\in\Z+\oh\,. 
\end{equation}
Consider the following function
\begin{equation}
  \label{qdil}
  \qd{z}=\prod_{n=0}^\infty (1-q^n z) \,.
\end{equation}
For various reasons, in particular because of the relation 
\eqref{qda}, this function is sometimes called the quantum
dilogarithm function \cite{FK}. 

It is clear that the function \eqref{Phip} has the
following expression in terms of the quantum dilogarithm 
$$
\Phi_{3D}(t,z)=
\begin{cases}
\dfrac{\qd{q^{1/2}/z}} 
{\qd{q^{1/2+t} z}}\,, & t \ge 0\,,\vspace{5pt} \\ 
\dfrac{\qd{q^{1/2-t}/z}}{\qd{q^{1/2} z}}\,, & t \le 0\,.
 \end{cases}
$$

\subsubsection{}\label{rtil}

\begin{figure}[htbp]\psset{unit=0.5cm}
  \begin{center}
    \begin{pspicture}(-7,-4)(7,5)
\scriptsize
\showgrid
\psset{dimen=middle}
\psaxes[axesstyle=frame,Ox=-7,Oy=-4,Dx=2,Dy=2,ticks=none](-7,-4)(7,5)
\psdiamond[fillstyle=solid,fillcolor=lightgray](0,4.5)(1,.5)
\psdiamond[fillstyle=solid,fillcolor=lightgray](1,2.0)(1,.5)
\psdiamond[fillstyle=solid,fillcolor=lightgray](2, .5)(1,.5)
\psdiamond[fillstyle=solid,fillcolor=lightgray](3,-1.0)(1,.5)
\psdiamond[fillstyle=solid,fillcolor=lightgray](4,-2.5)(1,.5)
\psdiamond[fillstyle=solid,fillcolor=lightgray](5,-3.0)(1,.5)
\psdiamond[fillstyle=solid,fillcolor=lightgray](6,-3.5)(1,.5)
\psdiamond[fillstyle=solid,fillcolor=lightgray](-1,3.0)(1,.5)
\psdiamond[fillstyle=solid,fillcolor=lightgray](0,1.5)(1,.5)
\psdiamond[fillstyle=solid,fillcolor=lightgray](1,-1.0)(1,.5)
\psdiamond[fillstyle=solid,fillcolor=lightgray](2,-1.5)(1,.5)
\psdiamond[fillstyle=solid,fillcolor=lightgray](3,-3.0)(1,.5)
\psdiamond[fillstyle=solid,fillcolor=lightgray](4,-3.5)(1,.5)
\psdiamond[fillstyle=solid,fillcolor=lightgray](-2,1.5)(1,.5)
\psdiamond[fillstyle=solid,fillcolor=lightgray](-1,0.0)(1,.5)
\psdiamond[fillstyle=solid,fillcolor=lightgray](0,-1.5)(1,.5)
\psdiamond[fillstyle=solid,fillcolor=lightgray](1,-3.0)(1,.5)
\psdiamond[fillstyle=solid,fillcolor=lightgray](2,-3.5)(1,.5)
\psdiamond[fillstyle=solid,fillcolor=lightgray](-3,0.0)(1,.5)
\psdiamond[fillstyle=solid,fillcolor=lightgray](-2,-1.5)(1,.5)
\psdiamond[fillstyle=solid,fillcolor=lightgray](-1,-3.0)(1,.5)
\psdiamond[fillstyle=solid,fillcolor=lightgray](0,-3.5)(1,.5)
\psdiamond[fillstyle=solid,fillcolor=lightgray](-4,-2.5)(1,.5)
\psdiamond[fillstyle=solid,fillcolor=lightgray](-3,-3.0)(1,.5)
\psdiamond[fillstyle=solid,fillcolor=lightgray](-2,-3.5)(1,.5)
\psdiamond[fillstyle=solid,fillcolor=lightgray](-5,-3.0)(1,.5)
\psdiamond[fillstyle=solid,fillcolor=lightgray](-4,-3.5)(1,.5)
\psdiamond[fillstyle=solid,fillcolor=lightgray](-6,-3.5)(1,.5)
\end{pspicture}
    \caption{Horizontal tiles of the tiling corresponding to diagram in Figure \ref{f1}}
    \label{f3}
  \end{center}
\end{figure}
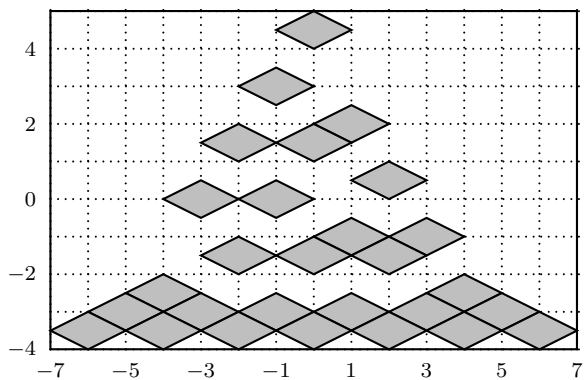

In order to make a better connection with the geometry of
3D diagrams, let us introduce a different encoding of
diagrams by subsets in the plane. Given a plane 
partition
$$
\pi=(\pi_{ij})\,, \quad i,j=1,2,\dots\,,
$$
we set
$$
\tfS(\pi)=\{(j-i,\pi_{ij}-(i+j-1)/2)\}\,,
\quad i,j=1,2,\dots \,.
$$
There is a well-known correspondence between 3D diagrams and
tilings of the plane by rhombi. Namely, the tiles are the 
images of the faces of the 3D diagram under the projection
\begin{equation}
(x,y,z)\mapsto (t,h)=(y-x,z-(x+y)/2)\,.\label{xyz}
\end{equation}
The tiling corresponding to the diagram in Figure \ref{f1}
is shown in Figure \ref{f3}. 

It is clear that under this 
correspondence the horizontal faces of a 3D diagram are
mapped to the horizontal tiles and that the positions 
of the horizontal tiles uniquely determine the tiling and
the diagram $\pi$. The set 
$$
\tfS(\pi) \subset \Z \times \oh \, \Z
$$
is precisely the set of the centers of the horizontal tiles. 
It is also clear that if $\cl$ corresponds to the diagram 
$\pi$, then
$$
(t,h) \in \tfS(\pi) \Leftrightarrow (t,x+|t|/2)\in \fS(\cl) \,.
$$
Theorem \ref{t1} specializes, therefore, to the following
statement

\begin{corollary}\label{cor3D}
For any set $\{(t_i,h_i)\}$, we have
$$
\Prob\left(\{(t_i,h_i)\} \subset \tfS(\pi)\right) =
\det\big[K_{3D} ((t_i,h_i),(t_j,h_j)) \big]\,,
$$
where the kernel $K_{3D}$ is given by the following formula
\begin{multline}
  \label{K3D}
  K_{3D} ((t_1,h_1),(t_2,h_2)) = \\
\frac1{(2\pi i)^2}
\int_{|z|=1\pm\epsilon} \int_{|w|=1\mp\epsilon} 
\frac1{z-w} \frac{\Phi_{3D}(t_1,z)}{\Phi_{3D}(t_2,w)}
\frac{dz \, dw}{z^{h_1+\frac{|t_1|+1}2} w^{-h_2-\frac{|t_2|-1}2}}\,.
\end{multline}
Here the function $\Phi_{3D}(t,z)$ is defined by
\eqref{Phip}, $0<\epsilon\ll 1$, and one picks the plus sign if $t_1 \ge t_2$
and the negative sign otherwise.
\end{corollary}

\section{Asymptotics}\label{secA}

\subsection{The local shape of a large 3D diagram}

\subsubsection{}

Our goal in this section is to illustrate how suitable 
is the formula \eqref{bK2} for asymptotic analysis. In 
order to be specific, we work out one concrete
example, namely the local shape of a 3D diagram
distributed according to the measure $\fM_q$ as $q\to 1$.
Our computations, however, will be of a very abstract and general 
nature and applicable to a much wider variety of
specialization.

The reader will notice that the passage to the 
asymptotics in \eqref{bK2} is so straightforward that even
the saddle-point analysis is needed in only a very weak 
form, namely as the statement that 
\begin{equation}\label{saddle}
\int_\gamma e^{M S(x)} \, dx  \to 0 \,, \quad M\to +\infty\,,
\end{equation}
provided the function $S(x)$ is smooth and $\Re S(x)< 0$ 
for all but finitely many points $x\in\gamma$.

\subsubsection{}
Let $q=e^{-r}$ and $r\to +0$. We begin with the following 

\begin{lemma} 
We have the following convergence in probability
$$
r^3 |\pi| \to
2\zeta(3) \,, \quad r\to +0 \,,
$$
where $|\pi|$ is 
the volume of a 3D diagram 
$\pi$ sampled from  the measure $\fM_q$. 
\end{lemma}
\begin{proof}
First consider the expectation of $|\pi|$ 
\begin{equation*}
\bE|\pi| = \frac{q\frac{d}{dq} Z_{3D}}{Z_{3D}} =
\sum_{n\ge 1} \frac{n^2 q^n}{1-q^n} = 
\sum_{n,k\ge 1} n^2 q^{nk} = 
\sum_k \frac{q^k(1+q^k)}{(1-q^k)^3} 
\sim
\frac{2\zeta(3)}{r^3}  \,. 
\end{equation*}
Similarly, the variance of $|\pi|$ behaves like
$$
\Var |\pi| = q\frac{d}{dq} \, \bE|\pi| =
o(r^{-6})\,,
$$
whence $\Var (r^3|\pi|) \to 0$, which concludes the proof. 
\end{proof}

\subsubsection{}

It follows that as $r\to +0$ the typical 3D diagram $\pi$,
scaled by $r$ in all directions, approaches the suitably scaled
limit shape for typical 3D diagrams of a large volume described in \cite{CK}. 
Below we will
also see this limit shape appear from our calculations. 

We are interested in the $r\to+0$ limiting local structure of $\pi$
in the neighborhood of various points in the limit shape. In 
other words, we are interested in the limit of the 
kernel \eqref{K3D} as 
$$
r t_i \to \tau \,, \quad r h_i \to \chi \,,
$$
where the variables 
$\tau$ and $\chi$ describe the global position on the limit
shape, 
in such a way that the relative distances
$$
\Dt=t_1-t_2\,, \quad \Dh=h_1-h_2  \,,
$$
remain fixed. This limit is easy to obtain by a combination 
of residue calculus with saddle-point argument. Since the 
measure $\fM_q$ is obviously symmetric with respect to
the reflection $t\mapsto -t$, we can without loss of generality
assume that $\tau \ge 0$ in our computations. 

\subsubsection{}

We have the following $r\to+0$ asymptotics:
\begin{equation}
  \label{qda}
 \ln \qd{z} \sim  r^{-1} \, \int_0^z \frac{\ln(1-w)}{w} \, dw =
- r^{-1} \, \dil(1-z) \,, 
\end{equation}
where 
$$
\dil(1-z)=\sum_{n} \frac{z^n}{n^2}\,, \quad |z|\le 1\,,
$$
analytically continued 
with a cut along $(1,+\infty)$. 

Introduce the following
function
$$
S(z;\tau,\chi)=-(\tau/2+\chi)\ln z - \dil(1-1/z) + 
\dil(1-e^{-\tau} z) 
$$
and recall that we made the assumption that $\tau\ge 0$. 
The function $S(z;\tau,\chi)$ is analytic in the complex
plane with cuts along $(0,1)$ and $(e^\tau,+\infty)$.

As $r\to +0$, the exponentially large term in the integrand
in \eqref{K3D} is
\begin{equation}
  \label{expb}
 \exp\left(\frac1r(S(z;\tau,\chi)-S(w;\tau,\chi))\right)\,.  
\end{equation}
The saddle-point method suggests, therefore, to look at the 
critical points of the function $S(z;\tau,\chi)$. Since
$$
z \frac{d}{dz} S(z;\tau,\chi) = -\tau/2 - \chi -
\log(1-1/z)(1-e^{-\tau} z) \,,
$$
the critical points of $S$ are the roots of the quadratic 
polynomial 
$$
(1-1/z)(1-e^{-\tau} z)=e^{-\tau/2-\chi}\,.
$$
The two roots of this polynomial are complex conjugate
if 
\begin{equation}
  \label{ineqchi}
|e^{\tau/2}+e^{-\tau/2}-e^{-\chi}| < 2\,,   
\end{equation}
which can be expressed equivalently as
\begin{equation}
- 2 \ln \left(2 \cosh \frac\tau4\right) <  
\chi < - 2 \ln \left(2 \sinh \frac\tau4\right) \,.
\label{ineqchi2}
 \end{equation}
In the case when the roots are complex conjugate they lie on the circle 
$$
\gamma=\{|z|=e^{\tau/2}\}
$$
As we shall see below, the inequality \eqref{ineqchi}
describes precisely the possible values 
of $(\tau,\chi)$ that correspond to the bulk of
the limit shape.

\subsubsection{}
The following elementary properties of the function 
$S(z;\tau,\chi)$
\begin{gather*}
  S(\bar z; \tau,\chi) = \overline{S(z;\tau,\chi)}\,,\\
  S(z;\tau,\chi) + S(e^{\tau}/z;\tau,\chi) = - (\tau/2+\chi) \tau \,,
\end{gather*}
imply that on the circle $\gamma$ the real part 
of $S$ is constant, namely,
$$
\Re S(z;\tau,\chi) = - (\tau/2+\chi) \tau/2 \,, \quad z\in\gamma\,.
$$
On $\gamma$ we also have
$$
z \frac{d}{dz} S(z;\tau,\chi) = 3\tau/2 - \chi - \ln|e^\tau-z|^2 \,,
\quad z\in\gamma\,,
$$

\begin{figure}[htbp]
\begin{center}
\begin{pspicture}(-3.000000,-4.000000)(5.000000,4.000000)
\uput[l](-2.117000,0){$\gamma$}
\uput{10pt}[u](.995800,1.868000){$z_c$}
\uput{10pt}[d](.995800,-1.868000){$\bar z_c$}
\uput[l](2.117000,0){$e^{\tau/2}$}
\uput[r](4.482000,0){$e^{\tau}$}
\pcline{<-*}(2.504500,3.426000)(4.482000,0)
\Aput{$e^{\frac34\tau-\frac12\chi}$}
\SpecialCoor\pscircle(0,0){2.117000}
\psarc(4.482000,0){3.955000}{120}{240}
\psline{->}(2.117000;0.000000)(3.146000;0.000000)
\psline{->}(2.117000;18.000000)(2.992000;18.000000)
\psline{->}(2.117000;36.000000)(2.646000;36.000000)
\psline{->}(2.117000;54.000000)(2.270000;54.000000)
\psline{->}(2.117000;72.000000)(1.938000;72.000000)
\psline{->}(2.117000;90.000000)(1.665000;90.000000)
\psline{->}(2.117000;108.000000)(1.452000;108.000000)
\psline{->}(2.117000;126.000000)(1.291000;126.000000)
\psline{->}(2.117000;144.000000)(1.180000;144.000000)
\psline{->}(2.117000;162.000000)(1.115000;162.000000)
\psline{->}(2.117000;180.000000)(1.093000;180.000000)
\psline{->}(2.117000;198.000000)(1.115000;198.000000)
\psline{->}(2.117000;216.000000)(1.180000;216.000000)
\psline{->}(2.117000;234.000000)(1.291000;234.000000)
\psline{->}(2.117000;252.000000)(1.452000;252.000000)
\psline{->}(2.117000;270.000000)(1.665000;270.000000)
\psline{->}(2.117000;288.000000)(1.938000;288.000000)
\psline{->}(2.117000;306.000000)(2.270000;306.000000)
\psline{->}(2.117000;324.000000)(2.646000;324.000000)
\psline{->}(2.117000;342.000000)(2.992000;342.000000)
\end{pspicture}
\caption{Gradient of $\Re S(z)$ on the circle $\gamma=\{|z|=e^{\tau/2}\}$}
\label{f4}
\end{center} \end{figure}
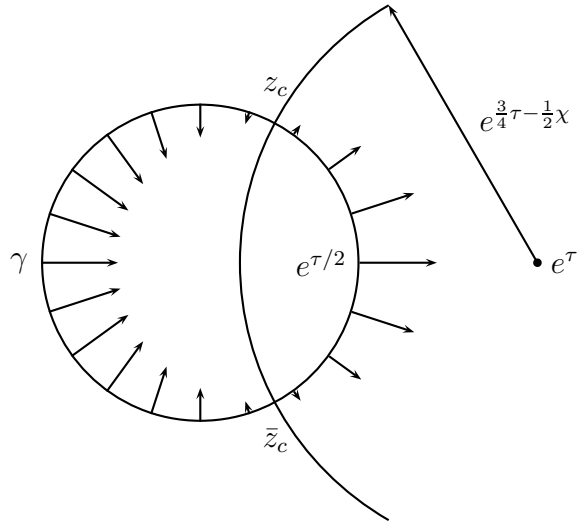

{}From this it is clear that when the critical points of $S$
are complex conjugate, they are the points of intersection of 
two following circles
\begin{equation}
  \label{intc}
  \{z_c,\bar z_c\}=\{|z|=e^{\tau/2}\} \cap 
\{|z-e^\tau|=e^{3\tau/4-\chi/2}\} \,.
\end{equation}
This is illustrated in Figure \ref{f4} which also shows the vector
field 
$$
\nabla\left(\Re S(z;\tau,\chi)\right) = \frac{z^2}{e^\tau} \, \frac{d}{dz} S(z) \,,
\quad z\in\gamma\,.
$$

\subsubsection{}

Now we are prepared to do the asymptotics in \eqref{K3D}. 
We can deform the contours of the integration in \eqref{K3D}
as follows
$$
 K_{3D} ((t_i,h_i),(t_j,h_j)) = 
\frac1{(2\pi i)^2}
\int_{(1\pm\epsilon)\gamma} dz  \int_{(1\mp\epsilon)\gamma} dw \quad
\dots
$$
where dots stand for the same integrand as in \eqref{K3D}, 
$0<\epsilon\ll 1$, and we pick the plus sign if $t_1 \ge t_2$
and the negative sign otherwise. 

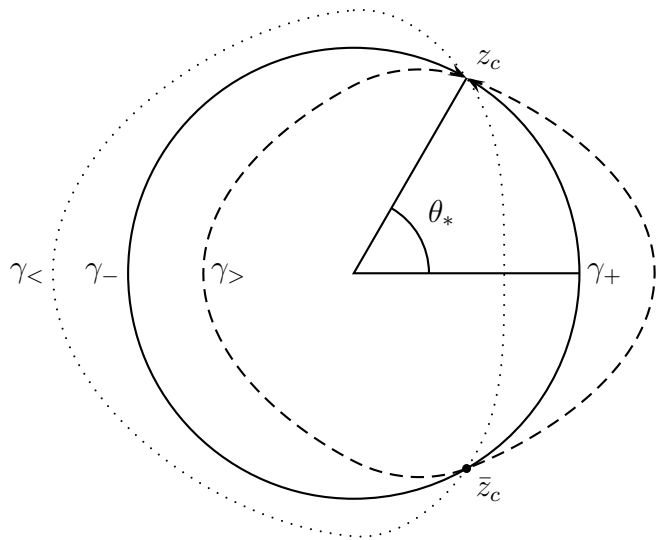
\begin{figure}[htbp]
  \begin{center}
\begin{pspicture}(-4,-4)(4,4)
\psarc[arrowsize=5pt]{<-*}(0,0){3}{60}{300}
\psarc[arrowsize=5pt]{*->}(0,0){3}{-60}{60}
\psccurve[linestyle=dashed](-2,0)(0,2.5)(1.5,2.5981)(4,0)(1.5,-2.5981)(0,-2.5)
\psccurve[linestyle=dotted]
(-4,0)(0,3.5)(1.5,2.5981)(2,0)(1.5,-2.5981)(0,-3.5)
\uput{3pt}[l](-4,0){$\gamma_<$}
\uput{3pt}[r](-2,0){$\gamma_>$}
\uput{3pt}[l](-3,0){$\gamma_-$}
\uput{3pt}[r](3,0){$\gamma_+$}
\uput{5pt}[ur](1.5,2.5981){$z_c$}
\uput{5pt}[dr](1.5,-2.5981){$\bar z_c$}
\psline(3,0)(0,0)(1.5,2.5981)
\psarc(0,0){1}{0}{60}
\SpecialCoor
\uput[ur](1;30){$\theta_*$}
\end{pspicture}

    \caption{Contours $\gamma_>$(dashed), $\gamma_<$(dotted), and
$\gamma_\pm$} \label{f5}
  \end{center}
\end{figure}

Now we define the contours $\gamma_>$, $\gamma_<$, $\gamma_+$,
$\gamma_-$. This definition will be illustrated by Figure \ref{f5}.
The contour $\gamma_>$ is the circle $|z|=e^{\tau/2}$ slightly
deformed in the direction of the gradient of $\Re S$, see Figure
\ref{f4}. Similarly, the contour $\gamma_{<}$ is the
same circle $|z|=e^{\tau/2}$ slightly pushed in the opposite 
direction. The contours $\gamma_\pm$ are the arcs of the 
circle $|z|=e^{\tau/2}$ between $\bar z_c$  and $z_c$, oriented
toward $z_c$. 

Deforming the contours, and picking the residue at $z=w$, we 
obtain
$$
K_{3D} ((t_i,h_i),(t_j,h_j)) = \int^{(1)} + \int^{(2)} \,,
$$
where
$$
\int^{(1)} = 
\frac1{(2\pi i)^2}
\int_{\gamma_<} dz  \int_{\gamma_>} dw \quad
\dots
$$
with the same integrand in as in \eqref{K3D} and
$$
\int^{(2)} = \frac1{2\pi i} \int_{\gamma_\pm} 
\frac
{\qd{q^{1/2+t_2} w}}
{\qd{q^{1/2+t_1} w}}
\frac{dw}{w^{\Dh+\Dt/2+1}} \,,
$$
where we choose $\gamma_+$ if $t_1  \ge t_2$ and 
$\gamma_-$ otherwise. As we shall see momentarily,
$$
\int^{(1)}  \to  0 \,, \quad r\to +0\,,
$$
while $\int^{(2)}$ has a simple limit. 

\subsubsection{}

It is obvious that 
$$
\int^{(2)} \to  \frac1{2\pi i} \int_{\gamma_\pm} 
(1-e^{-\tau} w)^{\Dt} 
\frac{dw}{w^{\Dh+\Dt/2+1}}\,, \quad r\to +0 \,.
$$
The change of variables $w=e^\tau w'$ makes this integral
a standard incomplete beta function integral:
$$
\int^{(2)} \to \frac{e^{-\tau(\Dh+\Dt/2)}}
{2\pi i} \int_{e^{-\tau} \gamma_\pm}
(1-w)^{\Dt} \frac{dw}{w^{\Dh+\Dt/2+1}}\,, \quad r\to +0 \,.
$$
Now notice that the prefactor $e^{-\tau(\Dh+\Dt/2)}$ will cancel 
out of any determinant with this kernel, so it can be ignored. 
We, therefore, make the following

\begin{definition}
  Introduce the following incomplete beta function kernel
  $$
  \bB_\pm(k,l;z)= \frac{1}{2\pi i}\int_{\bar z}^{z} (1-w)^{k}
  w^{-l-1} \, dw\,,
  $$
  where the path of the integration crosses $(0,1)$ for the plus
  sign and $(-\infty,0)$ for the minus sign.
\end{definition}

\subsubsection{}

It is also obvious from our construction that the function 
$\Re S(z)$ reaches its maximal value on the contour $\gamma_<$
precisely at the points $z_c$ and $\bar z_c$. The same 
points are the minima of the function $S(z)$ on the contour
$\gamma_>$. The behavior of the integral $\int^{(1)}$
is thus determined by the term \eqref{expb} and, by the 
basic principle \eqref{saddle}, the limit of $\int^{(1)}$ 
vanishes.

\subsubsection{}

We can summarize our discussion as follows. Set
$z_*=e^{-\tau} z_c$, in other words,
\begin{equation}
\label{zst}
  \{z_*, \bar z_*\} = \{|z|=e^{-\tau/2}\} \cap 
\{|z-1|=e^{-\tau/4-\chi/2}\}\,, \quad \Im z_*>0\,.
\end{equation}
It is  convenient to extend the meaning of 
$z_*$ to denote the point on the circle $|z|=e^{-\tau/2}$
which is the closest point to the circle $|z-1|=e^{-\tau/4-\chi/2}$
in case when the two circles do not intersect. In other
words, we complement the definition \eqref{zst} by  setting 
\begin{equation}
  \label{degz}
 z_*=
\begin{cases}
e^{-\tau/2} \,, & \quad  e^{-\chi/2} < e^{\tau/4}-e^{-\tau/4} \,, \\
-e^{-\tau/2}\,, & \quad  e^{-\chi/2} > e^{\tau/4}+e^{-\tau/4} \,. \\
\end{cases} 
\end{equation}

We have established the following

\begin{theorem}\label{t2}
Let $U=\{(t_i,h_i)\}$ and suppose that as $r\to+0$
$$
r t_i\to \tau\ge 0\,, \quad r h_i \to \chi\,,\quad i=1,2,\dots\,, 
$$
in such a way that the differences
$$
\Dt_{ij}=t_i-t_j\,, \quad \Dh_{ij}=h_i-h_j\,,
$$
remain fixed. Then, as $r\to+0$,
$$
\Prob\{U\subset\tfS(\pi)\}\to  
\det\left[
\bB_\pm\left(\Dt_{ij},\Dh_{ij}+\frac{\Dt_{ij}}2;z_*\right)
\right]\,,
$$
where the point $z_*=z_*(\tau,\chi)$ is defined in \eqref{zst}and \eqref{degz}
and the 
choice of the plus sign corresponds to $\Dt_{ij} =t_i - t_j\ge 0$. 
\end{theorem}

\begin{remark}
These formulas can be transformed, see Section \ref{iint}, into a
double integral of form 
considered in \cite{CKP}, Proposition 8.5 and Conjecture 13.5.
 \end{remark}

\begin{remark}
 Observe that the limit correlation are trivial in the cases
covered by \eqref{degz}. In other words, they are nontrivial
unless the inequality \eqref{ineqchi} is  satisfied. 
This means that the inequality \eqref{ineqchi} describes
the values of $(\tau,\chi)$ that correspond to the bulk of
the limit shape. 
 \end{remark}

\subsubsection{}

In particular, denote by $\rho_*(\tau,\chi)$ the
limit of $1$-point correlation function
$$
K_{3D}((t,h),(t,h)) \to \rho_*(\tau,\chi)\,, \quad 
rt\to\tau, \,  r h \to \chi\,. 
$$
This is the limiting density of the horizontal tiles at
the point $(\tau,\chi)$. We have the following 

\begin{corollary}\label{c2}
The limiting density of horizontal tiles is 
$$
\rho_*(\tau,\chi) = \frac{\theta_*}{\pi}  \,,
$$
where $\theta_*=\arg z_*$ (see Figure \ref{f5}), that is, 
\begin{equation}
\theta_*=
\arccos \left( \cosh\frac\tau2 - \frac{e^{-\chi}}2\right) \,. \label{th*}
\end{equation}
\end{corollary}

\begin{figure}[htbp]
  \begin{center}
    \scalebox{0.6}{\includegraphics{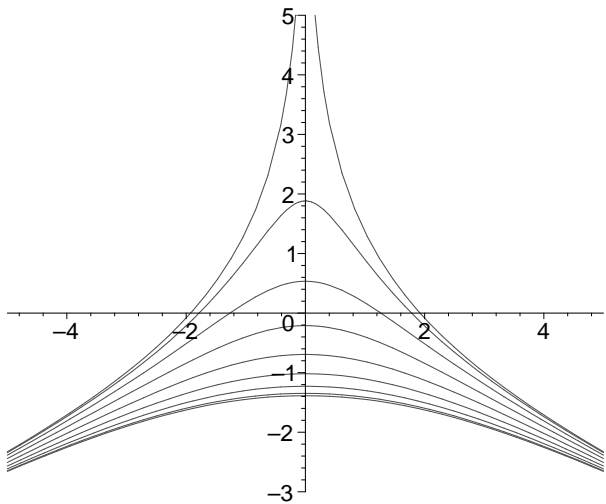}} 
    \caption{Level sets of the density of horizontal tiles}
    \label{f6}
  \end{center}
\end{figure}

The level sets of the density as functions of 
$\tau$ and $\chi$ are plotted in Figure \ref{f6}. More
precisely, Figure \ref{f6} shows the curves  
$$
\theta_*(\tau,\chi)= \frac{k\pi}8\,, \quad k=0,\dots,8\,.
$$
The knowledge of density $\rho_*$ is equivalent to the
knowledge of the limit shape, see Sections \ref{limsh1} and
\ref{limsh2}. We point out, however, that additional
analysis is needed to prove the convergence to the limit
shape in a suitable metric as in \cite{CK}.

\subsubsection{}

Corollary \ref{c2} can be generalized as follows: 

\begin{corollary}
The equal time correlations are given by the
discrete sine kernel, that is, if 
$t_1=t_2=\dots$ then 
$$
\Prob\{U\subset\tfS(\pi)\}\to  
\det\left[
\frac{\sin(\theta_*(h_i-h_j))}{\pi(h_i-h_j)}
\right]\,, 
$$
as $r\to+0$.
\end{corollary}

\subsubsection{}\label{iint}

The incomplete beta kernel $\bB_\pm$ can be transformed 
into a double integral of the form considered in \cite{CKP}
as follows. 

Using the following standard integral 
$$
\frac1{2\pi i}\int_{|z|=\alpha} \frac{z^{-k-1} \, dz}{1-\beta z} 
=
\begin{cases}
  \beta ^k \,, & k\ge 0,\, |\alpha \beta|<1 \,,\\
  -\beta^k \,, & k<0, \, |\alpha\beta|>1 \,, \\
  0\,, & \text{otherwise}\,, 
\end{cases}
$$
we can replace the condition 
$$
|w-1|\lessgtr e^{-\tau/4-\chi/2}\,,
$$
in the definition of $\bB_\pm(k,l;z_*)$ by an extra
integral. We obtain
\begin{equation}
\bB_\pm(k,l;z_*)=\frac1{(2\pi i)^2}
\iint_
{
  \begin{subarray}{l}
    |w|=e^{-\tau/2}\\ |z|=e^{\tau/4+\chi/2}
\end{subarray}
}
\frac{z^{-k-1} w^{-l-1}}{1-z+zw}\, dz\, dw\,,
\label{eii}
\end{equation}
where the plus sign corresponds to $k\ge 0$. 

\subsubsection{}\label{limsh1}

Let $z(\tau,\chi)$ denote the $z$-coordinate of 
the point on the limit shape corresponding to the point 
$(\tau,\chi)$. This function can be obtained
by integrating the density $\rho_*(\tau,\chi)$ as
follows. 

Consider a tiling such as the one in the Figure 
\ref{f3} and the corresponding 3D diagram, which
for the tiling in  Figure \ref{f3} is shown in 
Figure \ref{f1}. It is clear that the $z$-coordinate
of the face corresponding to a given horizontal tile equals 
the number of holes (that is, positions not occupied by a horizontal
tile) below it.  It follows that
\begin{equation}
z(\tau,\chi)=\int_{-\infty}^\chi (1-\rho_*(\tau,s)) \, ds \,. \label{z}
\end{equation}
Here, of course, the lower limit of integration can 
be any number between $-\infty$ and the lower boundary of the
limit shape given by the equation \eqref{ineqchi2}. 

\begin{figure}[htbp]
  \begin{center}
   \scalebox{0.6}{\includegraphics{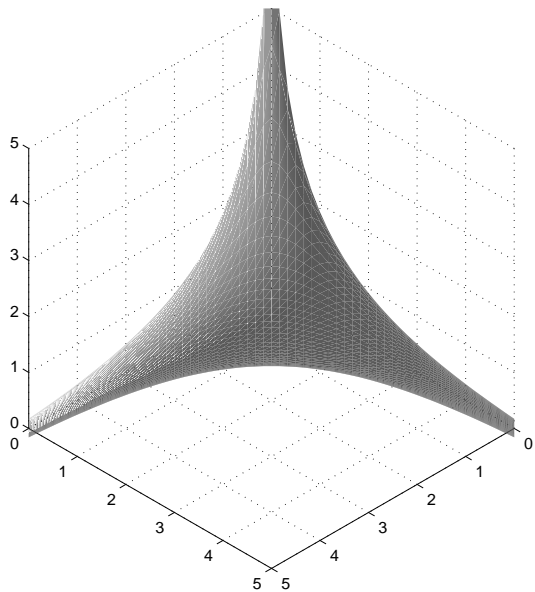}}     
    \caption{The limit shape}\label{f7}
     \end{center}
\end{figure}

Integrating \eqref{th*} we obtain the following formula
\begin{equation}
z(\tau,\chi)=\frac1\pi\int_0^{\pi-\theta_*}
\frac{s \, \sin s \,ds}{\cos s + \cosh \frac\tau2}\,.\label{z1}
\end{equation}
This integral can be evaluated in terms of the 
dilogarithm function. From \eqref{xyz} we can now compute 
the other two coordinates as follows
\begin{equation}
x(\tau,\chi)=z(\tau,\chi)-\chi-\frac\tau2\,,
\quad
y(\tau,\chi)=z(\tau,\chi)-\chi+\frac\tau2\,,\label{xy}
\end{equation}
which gives a parametrization of the limit shape. 
A plot of the limit shape is shown in Figure \ref{f7}

\subsubsection{}\label{limsh2}

To make the connection to the parametrization of the 
limit shape given in \cite{CK}, 
let us now substitute for $\rho_*$ in the integral
\eqref{z} the formula \eqref{eii}.
After a simple coordinate change we obtain 
$$
\rho_*(\tau,\chi) = \frac1{4\pi^2} \iint_0^{2\pi}
\frac{du \,dv}{1+e^{\chi/2+\tau/4+iu}+e^{\chi/2-\tau/4+iv}} \,.
$$
Integrating this in $\chi$ we obtain
$$
z(\tau,\chi)=\frac1{2\pi^2} \iint_0^{2\pi} 
\ln \left|1+e^{\chi/2+\tau/4+iu}+e^{\chi/2-\tau/4+iv}\right| \,
du \,dv \,.
$$
This together with \eqref{xy} is equivalent to following parametrization 
of the limit shape found in \cite{CK}
$$
(x,y,z)=\left(f(A,B,C)-2\ln A,f(A,B,C)-2\ln B, f(A,B,C)-2 \ln C
\right)\,, 
$$
where $A,B,C>0$ and 
$$
f(A,B,C)=\frac1{2\pi^2} \iint_0^{2\pi} 
\ln \left|A+ B e^{iu}+C e^{iv}\right| \,
du \,dv \,.
$$
This parametrization is manifestly symmetric in $x$, $y$, and $z$.
However, it involves integrals that are more complicated to 
evaluate than the parametrization \eqref{z1}. 

Note that the shape considered in \cite{CK} differs from 
our by a factor of $2$ due to different scaling conventions.

\subsection{Universality}

\subsubsection{}
The reader has surely noticed that all that we really needed
for the asymptotics is to be able to deform the 
contours of integration as in Figure \ref{f5} so that
the points of the intersection of $\gamma_>$ and $\gamma_<$
the minima of the function $\Re S(z)$ on one curve and 
the maxima --- on the other. The intersection points
are then forced to be the critical points of the function
$S(z)$. 

It is clear that this 
principle is very general and applicable in a potentially 
very large variety of situations, such as, for example,
in the case of Plancherel measure, see below. In particular,
the existence of contours of a certain kind is a property
which is preserved under small perturbations. 

\subsubsection{}
One such perturbation would be to consider anisotropic 
partitions. The anisotropy in the $t$-direction is 
especially easy to introduce: one just should replace
$q^{|m|}$ in \eqref{specp} by $q^{V(m)}$ for some 
function $V$.

\subsubsection{}\label{eqtime}
Observe that the asymptotics of the equal time correlations
are especially easy to obtain in our approach. This is 
because  
$$
\Res_{z=w}\frac{dz}{z-w} \, \frac{\Phi(t_1,z)}
{\Phi(t_2,w)}=1\,, \quad t_1=t_2\,,
$$
and hence the analog of the integral $\int^{(2)}$ becomes simply
the integral
$$
\int^{(2)} = \frac1{2\pi i} \int_{\gamma_\pm} 
\frac{dw}{w^{\Delta x+1}} \,,
$$
which leads to the discrete sine kernel in the variable 
$\Delta x$. 

\subsubsection{}
Let us illustrate these general remarks by briefly discussing 
 the asymptotics for
the poissonized Plancherel measure \cite{BOO}. It corresponds
to the following specialization of the Schur process
$$
\phi_{\Pl}[m]=
\begin{cases}
e^{\sqrt{\alpha} \, z}\,, & m=-\frac12\,,\\
e^{\sqrt{\alpha} \, z^{-1}}\,, & m=\frac12\,,\\
1\,, & \text{otherwise}\,,
\end{cases}
$$
where $\alpha>0$ is the poissonization parameter. 
In particular, only the partition $\lambda(0)$ is 
nontrivial. The correlation kernel at $t=0$ specializes
to 
$$
K_{\Pl}(x,y)=\frac1{(2\pi i)^2} \iint 
\frac{z^{-x-1/2} w^{y-1/2}}{z-w} \, e^{\sqrt{\alpha} (z-z^{-1}-w+w^{-1}) } \,
dz\,dw \,,
$$
which can be expressed in terms of Bessel functions of the 
argument $2\sqrt{\alpha}$, see \cite{BOO,J1}.

\subsubsection{}

 The $\alpha\to\infty$ asymptotics
of the kernel $K_{\Pl}$ is easy to obtain from the 
classical asymptotics of the Bessel functions, see \cite{BOO}. 
It is, however, instructive to see 
how  this can be done even more quickly in our
framework. Assume that 
$$
\frac{x}{\sqrt{\alpha}}, \frac{y}{\sqrt{\alpha}} \to \xi
$$
in such a way that $\Delta=x-y$ remains fixed. The critical 
points of the action 
$$
S_{\Pl} = z-z^{-1}-\frac{\xi}{2}\, \ln z 
$$
are complex precisely when $|\xi|<2$,
in which case they are the points $e^{\pm i\theta}$, where 
$$
\theta = \arccos(\xi/2) \,.
$$
So, the same argument as we employed above immediately yields
the following formula from \cite{BOO}.
$$
K_{\Pl}(x,y) \to \frac{\sin \theta \Delta}{\pi \Delta}\,.
$$

\subsubsection{}

Of course, a finer analysis (which was carried out in \cite{BOO})
is needed to justify depoissonization in the asymptotics.
 In our situation, a similar problem is to pass from the 
$q\to 1$ asymptotics of the measure $\fM_q$ to the asymptotics
of the uniform measures on partitions of a given volume $N$
as $N\to\infty$. In other words, further work is needed to 
verify the equivalence of ensembles in the asymptotics. 

\subsubsection{}

Similarly, finer analysis is needed to work with the asymptotics
at the edges of the limit shapes where one expects to
see the Airy kernel appear. In the edge scaling, the following
equivalent version of the formula \eqref{bK2}
\begin{equation}\label{sumK}
K((t_1,x_1),(t_2,x_2))=
\sum_{m=1/2}^{\infty}
[z^{x_1+m} w^{-x_2-m}] \, \frac{\Phi(t_1,z)}
{\Phi(t_2,w)}\,,
\end{equation}
may be helpful because, by analogy to the situation with the
Plancherel measure \cite{BOO,J1}, one could expect to see this 
sum to become the integral representing the Airy kernel
$$
\frac{\Ai(x) \Ai'(y) - \Ai'(x) \Ai(y)}{x-y}=
\int_0^\infty \Ai(x+s) \, \Ai(y+s) \, ds \,.
$$
Further discussion of the Airy-type asymptotics of the
integrals of the form \eqref{bK2}, \eqref{K3D} can be
found in \cite{O2}.

\subsubsection{}

Recently, the techniques of this paper were used in 
\cite{FS} to prove that, indeed, the boundary of a random 3D
Young diagram converges to the Airy process. 

\appendix

\section{Summary of the infinite wedge formulas}

Let the space $V$ be spanned by $\ul{k}$, $k\in\Z+\oh$.
The space $\LV$  is, by definition, spanned by vectors 
$$
v_S=\ul{s_1} \wedge \ul{s_2} \wedge  \ul{s_3} \wedge  \dots\,,
$$
where $S=\{s_1>s_2>\dots\}\subset \Z+\oh$ is such a subset that
both sets
$$
S_+ = S \setminus \left(\Z_{\le 0} - \oh\right) \,, \quad
S_- = \left(\Z_{\le 0} - \oh\right) \setminus S 
$$
are finite. We equip $\LV$ with the inner product 
in which the basis $\{v_S\}$ is orthonormal. In particular,
we have the vectors
$$
v_\lambda=\ul{\lambda_1-\tfrac12} \wedge \ul{\lambda_2-\tfrac32}
\wedge \ul{\lambda_4-\tfrac52} \wedge \dots \,,
$$
where $\lambda$ is a partition. 
The vector 
$$
\vac = \ul{-\tfrac12} \wedge \ul{-\tfrac32} \wedge \ul{-\tfrac52} \wedge \dots
$$
is called the vacuum vector. 

The operator $\psi_k$ is the exterior multiplication by $\ul{k}$
$$
\psi_k \left(f\right) = \ul{k} \wedge f  \,. 
$$
The operator $\psi^*_k$ is the adjoint operator.
These operators satisfy the canonical anti-commutation relations
$$
\psi_k \psi^*_k + \psi^*_k \psi_k = 1\,,
$$
all other anticommutators being equal to $0$. We have 
\begin{equation}\label{yn}
\psi_k \psi^*_k \,\, v_S = 
\begin{cases}
v_S\,, & k \in S \,, \\
0 \,, & k \notin S \,.
\end{cases}
\end{equation}

The operators $\al_n$ defined by 
$$
\al_n = \sum_{k\in\Z+\frac12} \psi_{k-n}\, \psi^*_k \,, \quad n=\pm 1,\pm2, \dots \,, 
$$
satisfy the Heisenberg commutation relations
$$
\left[\al_n, \al_m\right] = n \, \delta_{n,-m} \,,
$$
see the formula. Clearly, $\al^*_n=\al_{-n}$. It is clear from definitions that %
\begin{equation}\label{adal}
[\al_n,\psi(z)]=z^n\, \psi(z)\,, \quad  [\al_n,\psi^*(w)]=- w^{n}\, \psi^*(w) 
\end{equation}
and also that 
\begin{equation}
  \label{ann}
 \al_n \, \vac  = 0 \,, \quad  n\le 0 \,.  
\end{equation}

\end{document}